# Isotone Cones in Banach Spaces and Applications to Best Approximations of Operators without Continuity Conditions

Jinlu Li
Department of Mathematics
Shawnee State University
Portsmouth, Ohio 45662
USA

**Abstract**. In this paper, we introduce the concept of isotone cones in Banach spaces. Then we apply the order monotonic property of the metric projection operator to prove the existence of best approximations for some operators without continuity conditions in partially ordered Banach spaces.

## 1. Introduction

In linear and nonlinear analysis on Banach spaces, which includes fixed point theory, optimization theory, approximation theory, equilibrium theory, and variational analysis, etc., the metric projection and duality mapping are two especially important operators. They have been widely used to: 1) proving the solvability of some problems in these fields (see [1−5], [13−17]); 2) finding some algorithms to estimating the solutions (see [7−9], [18]).

A given Banach space, in addition to the algebraic and norm (topology) structures, can be equipped with a partially ordered structure. Such a partial order structure on a Banach space can be naturally induced by a closed and convex cone in this space. The selections of cones should be based on the given problems (see [6]). In analysis of Banach spaces, there are significant advantages of introducing partial orders: 1) the underlying spaces can be extended, for example, compact subsets may be replaced by chain-complete subsets; 2) the continuity conditions on the considered mappings may be substituted by order-monotonic conditions. To explain these ideas, we list the following fixed point theorem as an example and that is a consequence of Theorems 3.10 in [12]:

**Theorem**. *Let $(X, \succcurlyeq)$ be a partially ordered reflexive Banach space, in which the partial order $\succcurlyeq$ is induced by a nonempty closed and convex cone K, and D a nonempty bounded closed convex subset of X. Let $F: D \to 2^D \setminus \{\varnothing\}$ be a set-valued mapping satisfying the following conditions*:

(a). *F is $\succcurlyeq$-increasing upward*;
(b). *$F(x)$ is a nonempty closed and convex subset of D, for every $x \in D$*;
(c). *There are elements $y \in D$ and $v \in F(y)$ with $y \preccurlyeq v$.*

*Then, F has a fixed point and the set of fixed points of F is an inductive poset.*

In this theorem, the underlying space $D$ is a bounded closed convex subset of $X$ that is not required to be compact. In the reflexive Banach space $X$, $D$ is weakly compact. It insures that $(D, \succcurlyeq)$ is an $\succcurlyeq$-chain-complete poset. The considered mapping $F$ in this theorem is required to hold some ordering properties A1-A3 and it is not required to satisfy any continuity condition.

Let $P_K$ be the metric projection of $X$ onto $K$. As we mentioned in the first paragraph, $P_K$ is very useful in analysis on Banach spaces. In particular, $P_K$ has been used by many authors to construct some iteration processes to approximating fixed points of some mappings. In case, after introducing a partial order $\succcurlyeq$ on the give Banach space, in order to apply the above theorem for $P_K$, one will immediately ask:

$$\text{Is } P_K \succcurlyeq\text{-increasing?} \tag{1}$$

In this paper, we will answer the question (1) by constructing a counter example in $R^3$ to show:

$$\text{In general, } P_K \text{ is not } \succcurlyeq\text{-increasing.}$$

On the other hand, in [9], Kong, Liu and Wu provide a special closed and convex cone $K$ in Hilbert spaces such that the metric projection operator $P_K$ is $\succcurlyeq$-increasing, where $\succcurlyeq$ is induced by this cone $K$.

In this paper, for answering the question (1), we introduce the concept of "isotone cones" in Banach spaces, for which the metric projection operator is order-increasing. By applying this concept, we prove some fixed point theorems. As applications, we study the best approximation problems in Banach spaces.

2. **Preliminaries**

Let $X$ be a Banach space and $X^*$ its dual. The duality mapping $J: X \to 2^{X^*}$ is defined as follows: for every $x \in X$

$$J(x) = \{j \in X^*: \langle x, j \rangle = \|x\|^2 = \|j\|^2\}.$$

Let $K$ be a nonempty closed and convex cone in $X$. Let $\succcurlyeq$ be the partial order on $X$ induced by $K$; that is,

$$x \succcurlyeq y \text{ if and only if } x - y \in K, \text{ for all } x, y \in X.$$

Then $(X, \succcurlyeq)$ is a partially ordered Banach space. $(X, \succcurlyeq)$ is said to be a Banach lattice, if it satisfies:

(a) $\succcurlyeq$ is a lattice ordering on $X$, that is, for every $x, y \in X$, both $x \vee y$ and $x \wedge y$ exist;
(b) $\succcurlyeq$ is comparable to the norm $\|\cdot\|$, that is, for every $x, y \in X$,

$$0 \preccurlyeq x \preccurlyeq y, \text{ implies } \|x\| \leq \|y\|.$$

**Remarks**: Some authors say that $(X, \succcurlyeq)$ is a normal partially ordered Banach space, if it has property (b). It follows that every Banach lattice is a normal partially ordered Banach space.

Let $P_K: X \to 2^K$ be the metric projection of $X$ onto $K$ that is defined by, for every $x \in X$

$$P_K(x) = \{x' \in K: \|x-x'\| = \inf\{\|x-z\|: z \in K\}\}.$$

By applying the properties of the metric projection $P_K$, we provide the connections between $P_K$ and the duality mapping $J$ as below:

**Lemma 2.1.** $x' \in P_K(x)$ if and only if, there is $j(x' - x) \in J(x' - x)$ such that

$$\langle z - x', j(x' - x) \rangle \geq 0, \text{ for all } z \in K.$$

*Proof.* This lemma immediately follows from Theorem 3.8.4 part (iv) in [19]. □

In general, a nonempty closed and convex cone $K$ in a real vector space $X$ determines a partial order $\succcurlyeq$ on $X$, with respect to which, $(X, \succcurlyeq)$ may or may not be a vector lattice. In 1962, Kendall proved the known Choquet-Kendall Theorem in [11], that provides the criteria for partial orders in real vector spaces induced by cones to be lattice ordering relations. We first recall some concepts in [11] used in the Choquet-Kendall Theorem.

Let $D$ be a nonempty convex subset in a real vector space $X$. We write

$$D^* = \{\alpha x + (1-\alpha)y: \alpha \text{ real}; x, y \in D\}.$$

$D^*$ is called the minimal affine extension of $D$. Let

$$K = \{\lambda x: \lambda \geq 0; x \in D\}.$$

If $0 \notin D^*$, then $K$ is a convex cone in $X$, that is said to be generated by a base $D$.

A nonempty subset $D$ in a real vector space $X$ is a simplex when the intersection of two positively homothetic images of $D$,

$$x + \lambda D, y + \mu D \quad (x, y \in X; \lambda, \mu \geq 0)$$

is empty, or is a set $a + \gamma D$ ($a \in X; \gamma \geq 0$) of the same nature.

**The Choquet-Kendall Theorem.** *Let $D$ be a non-empty convex set of a real vector space $X$, and suppose that its minimal affine extension does not contain the zero vector. Let $C$ be the pointed positive convex cone $=\{\lambda b: \lambda \geq 0, b \in D\}$ having $D$ as base (of $X$). Then $X$, partly ordered by $C$, will be a vector lattice if and only if both*

(1°) *$D$ is a simplex, and*
(2°) *each B-segment $\{\lambda a + (1-\lambda)b: 0 \leq \lambda \leq 1; a, b \in D\}$ is contained in a maximal D-segment.*

## 3. A counter example

In this section, by using Lemma 2.1, we provide a counterexample to demonstrate that in a partially ordered Hilbert space $(X, \succcurlyeq)$, in which the partial order $\succcurlyeq$ is induced by a nonempty closed and convex cone $K \subseteq X$, the metric projection operator $P_K: X \to K$ is not $\succcurlyeq$-increasing.

Let $R^3$ be the 3-d Euclidean space. For every nonzero vectors $u, v \in R^3$, let $\angle(u, v)$ denote the angle formed by $u$ and $v$ satisfying

$$0 \leq \angle(u, v) \leq \pi.$$

**Example 3.1.** In $R^3$, take the vector $n = (1, 0, 0)$ and define a subset $K \subseteq R^3$ as follows:

$$K = \{\theta\} \cup \{u \in R^3 : 0 \leq \angle(u, n) \leq \frac{\pi}{4}\}.$$

Then $K$ is a closed and convex cone that induces a partial order $\succcurlyeq$ on $R^3$ and the metric projection operator $P_K$ from $R^3$ onto $K$ is not $\succcurlyeq$-increasing.

*Proof.* First of all, we show that the partial order $\succcurlyeq$ induced by $K$ is not a lattice ordering. To this end, take

$$D = \{(1, y, z) \in R^3 : y^2 + z^2 \leq 1\}.$$

One can check that

(i) $D$ is a nonempty convex subset in $R^3$;
(ii) The minimal affine extension of $D$ does not contain 0;
(iii) The closed convex cone $K$ is generated by this base $D$.

Take two positively homothetic images of $D$,

$$D = (0, 0, 0) + D \quad \text{and} \quad (0, \tfrac{1}{2}, 0) + D.$$

The intersection of the above two positively homothetic images of $D$ is a vesica piscis, that is neither empty, nor a set $a + \gamma D$ ($a \in R^3$; $\gamma \geq 0$) of the same nature. So $D$ is not a simplex in $R^3$ and the condition (1°) in The Choquet-Kendall Theorem is not satisfied. So $\succcurlyeq$ is not a lattice ordering. Hence $(R^3, \succcurlyeq)$ is not a Hilbert lattice.

From the inner product property, the inequalities $0 \leq \angle(u, v) \leq \pi$ are equivalent to

$$1 \geq \frac{\langle u, v \rangle}{\|u\|\|v\|} \geq -1, \text{ for every nonzero vectors } u, v \in R^3.$$

And the inequalities $0 \leq \angle(u, n) \leq \frac{\pi}{4}\}$ are equivalent to

$$1 \geq \frac{\langle u, v \rangle}{\|u\|\|v\|} \geq \frac{\sqrt{2}}{2}, \text{ for every nonzero vectors } u, v \in R^3.$$

Noticing that $\frac{\langle u, v \rangle}{\|u\|\|v\|} \leq 1$, for every nonzero vectors $u, v \in R^3$, then the cone $K$ is equivalently defined by

$$K = \{\theta\} \cup \left\{v \in R^3 : \frac{\langle v, n \rangle}{\|v\|\|n\|} \geq \frac{\sqrt{2}}{2}\right\}.$$

From $n = (1, 0, 0)$, the above equation can be rewritten as

$$K = \{\theta\} \cup \left\{ v = (x, y, z) \in R^3 : \frac{x}{\sqrt{x^2 + y^2 + z^2}} \geq \frac{\sqrt{2}}{2} \right\}.$$

It is easy to check that $(R^3, \succcurlyeq)$ is a partially ordered Hilbert space, that satisfies that, for every $u$, $v \in R^3$,

$$v \succcurlyeq u \quad \text{if and only if} \quad v - u \in K, \tag{2}$$

It is equivalent to

$$v \succcurlyeq u \quad \text{if and only if} \quad \frac{\langle v - u, n \rangle}{\|v - u\|} \geq \frac{\sqrt{2}}{2}. \tag{3}$$

The metric projection operator $P_K$ from $R^3$ onto $K$ clearly satisfies:

$$P_K(u) = \begin{cases} u, & \text{if } u \in K, \\ \theta, & \text{if } u \in -K. \end{cases}$$

Next we deduce a formula for $P_K$ on $R^3 \backslash (K \cup (-K))$. For any $u = (a, b, c) \in R^3$, we have that

$$u \in R^3 \backslash (K \cup (-K)) \quad \text{if and only if} \quad \frac{\pi}{4} < \angle(u, n) < \frac{3\pi}{4}.$$

The above inequalities are equivalent to

$$-\frac{\sqrt{2}}{2} < \frac{\langle u, n \rangle}{\|u\|} < \frac{\sqrt{2}}{2}, \text{ for any } u \in R^3 \backslash (K \cup (-K))$$

That is

$$-\frac{\sqrt{2}}{2} < \frac{a}{\sqrt{a^2 + b^2 + c^2}} < \frac{\sqrt{2}}{2}, \quad \text{for any } (a, b, c) \in R^3 \backslash (K \cup (-K)). \tag{4}$$

It must be satisfied that $\sqrt{b^2 + c^2} > 0$. Then it is induced to

$$(a, b, c) \in R^3 \backslash (K \cup (-K)) \quad \text{if and only if} \quad \sqrt{b^2 + c^2} > 0 \text{ and } -1 < \frac{a}{\sqrt{b^2 + c^2}} < 1. \tag{5}$$

For any given $u = (a, b, c) \in R^3 \backslash (K \cup (-K))$, let $w = P_K(u)$. The formula for $P_K(u)$ is

$$w = P_K(u) = P_K((a, b, c))$$

$$= \left( a + \left( \frac{1}{2} - \frac{a}{2\sqrt{b^2 + c^2}} \right) \sqrt{b^2 + c^2}, \left( \frac{a}{2\sqrt{b^2 + c^2}} + \frac{1}{2} \right) b, \left( \frac{a}{2\sqrt{b^2 + c^2}} + \frac{1}{2} \right) c \right)$$

$$= \left( \frac{a+\sqrt{b^2+c^2}}{2}, \frac{b(a+\sqrt{b^2+c^2})}{2\sqrt{b^2+c^2}}, \frac{c(a+\sqrt{b^2+c^2})}{2\sqrt{b^2+c^2}} \right)$$

$$= \frac{a+\sqrt{b^2+c^2}}{2}\left(1, \frac{b}{\sqrt{b^2+c^2}}, \frac{c}{\sqrt{b^2+c^2}}\right).$$

By equivalently verifying

$$\frac{\langle w, n \rangle}{\|w\|} = \frac{\sqrt{2}}{2},$$

it can be checked that

$$\angle(w, n) = \frac{\pi}{4}.$$

It implies that $w \in K$. Since $R^3$ is a Hilbert space, to show that $w = P_K(u)$, from Lemma 2.1, it is sufficient to verify that

$$\langle v - w, w - u \rangle \geq 0, \text{ for all } v \in K. \tag{6}$$

At first, we calculate

$$\langle w, w - u \rangle$$

$$= \left\langle \frac{a+\sqrt{b^2+c^2}}{2}\left(1, \frac{b}{\sqrt{b^2+c^2}}, \frac{c}{\sqrt{b^2+c^2}}\right), \frac{a+\sqrt{b^2+c^2}}{2}\left(1, \frac{b}{\sqrt{b^2+c^2}}, \frac{c}{\sqrt{b^2+c^2}}\right) - (a,b,c) \right\rangle$$

$$= \frac{a+\sqrt{b^2+c^2}}{2}\left\langle \left(1, \frac{b}{\sqrt{b^2+c^2}}, \frac{c}{\sqrt{b^2+c^2}}\right), \frac{a+\sqrt{b^2+c^2}}{2}\left(1, \frac{b}{\sqrt{b^2+c^2}}, \frac{c}{\sqrt{b^2+c^2}}\right) \right\rangle$$

$$- \frac{a+\sqrt{b^2+c^2}}{2}\left\langle \left(1, \frac{b}{\sqrt{b^2+c^2}}, \frac{c}{\sqrt{b^2+c^2}}\right), (a,b,c) \right\rangle$$

$$= 0. \tag{7}$$

Hence (6) becomes

$$\langle v, w - u \rangle \geq 0, \text{ for all } v \in K. \tag{8}$$

Let $v = (x, y, z) \in K$. From the definition of $K$, we have

$$1 \geq \frac{x}{\sqrt{x^2+y^2+z^2}} \geq \frac{\sqrt{2}}{2}.$$

It implies

$$x \geq \sqrt{y^2 + z^2}.$$

From (5), it implies

$$\frac{a}{\sqrt{b^2 + c^2}} - 1 < 0.$$

Then we have

$$\langle v, w - u \rangle$$

$$= \left\langle (x, y, z), \frac{a + \sqrt{b^2 + c^2}}{2}\left(1, \frac{b}{\sqrt{b^2 + c^2}}, \frac{c}{\sqrt{b^2 + c^2}}\right) - (a, b, c) \right\rangle$$

$$= \frac{\sqrt{b^2 + c^2}}{2}x + \frac{ab}{2\sqrt{b^2 + c^2}}y + \frac{ac}{2\sqrt{b^2 + c^2}}z - \frac{ax}{2} - \frac{by}{2} - \frac{cz}{2}$$

$$= \frac{1}{2}\left(\sqrt{b^2 + c^2}\,x + \frac{ab}{\sqrt{b^2 + c^2}}y + \frac{ac}{\sqrt{b^2 + c^2}}z - ax - by - cz\right)$$

$$= \frac{1}{2}\left(\left(\sqrt{b^2 + c^2} - a\right)x + \left(\frac{a}{\sqrt{b^2 + c^2}} - 1\right)(by + cz)\right)$$

$$\geq \frac{1}{2}\left(\left(\sqrt{b^2 + c^2} - a\right)x + \left(\frac{a}{\sqrt{b^2 + c^2}} - 1\right)\sqrt{b^2 + c^2}\sqrt{y^2 + z^2}\right)$$

$$\geq \frac{1}{2}\left(\left(\sqrt{b^2 + c^2} - a\right)\sqrt{y^2 + z^2} + \left(\frac{a}{\sqrt{b^2 + c^2}} - 1\right)\sqrt{b^2 + c^2}\sqrt{y^2 + z^2}\right)$$

$$= 0.$$

It proves (8). Then (6) follows from (7) and (8) immediately. Hence $w = P_K(u)$.

To summarize the above calculation, we obtained a formula for the projection operator $P_K$ from $R^3 \setminus (K \cup (-K))$ onto $K$: for any $(a, b, c) \in R^3 \setminus (K \cup (-K))$ satisfying

$$-\frac{\sqrt{2}}{2} < \frac{a}{\sqrt{a^2 + b^2 + c^2}} < \frac{\sqrt{2}}{2} \quad \text{or} \quad -1 < \frac{a}{\sqrt{b^2 + c^2}} < 1, \tag{5}$$

we have

$$P_K(a, b, c) = \left(\frac{a + \sqrt{b^2 + c^2}}{2}, \frac{b\left(a + \sqrt{b^2 + c^2}\right)}{2\sqrt{b^2 + c^2}}, \frac{c\left(a + \sqrt{b^2 + c^2}\right)}{2\sqrt{b^2 + c^2}}\right)$$

$$= \frac{a + \sqrt{b^2 + c^2}}{2}\left(1, \frac{b}{\sqrt{b^2 + c^2}}, \frac{c}{\sqrt{b^2 + c^2}}\right). \tag{9}$$

Next we provide some vectors in $R^3\backslash(K\cup(-K))$ to demonstrate that the metric projection operator $P_K$ is not $\succcurlyeq$-increasing. To this end, we take

$$u = \left(1, \sqrt{2}, 0\right) \quad \text{and} \quad v = \left(0, \frac{1}{\sqrt{2}}, \frac{1}{\sqrt{2}}\right).$$

It can be checked that both of $u$ and $v$ satisfy (4) or (5). Hence $u, v \in R^3\backslash(K\cup(-K))$. Calculate

$$\frac{\langle u-v, n\rangle}{\|u-v\|} = \frac{1}{\sqrt{1+\left(\frac{1}{\sqrt{2}}\right)^2 + \left(-\frac{1}{\sqrt{2}}\right)^2}} = \frac{\sqrt{2}}{2}.$$

From (2) or (3), it implies $u - v \in K$, that is,

$$u \succcurlyeq v. \tag{10}$$

From formula (9), we have

$$P_K(u) = \left(\frac{1+\sqrt{2}}{2}, \frac{1+\sqrt{2}}{2}, 0\right) \quad \text{and} \quad P_K(v) = \left(\frac{1}{2}, \frac{1}{2\sqrt{2}}, \frac{1}{2\sqrt{2}}\right).$$

It implies

$$\frac{\langle P_K(u) - P_K(v), n\rangle}{\|P_K(u) - P_K(v)\|}$$

$$= \frac{\frac{\sqrt{2}}{2}}{\sqrt{\left(\frac{\sqrt{2}}{2}\right)^2 + \left(\frac{1+\sqrt{2}}{2\sqrt{2}}\right)^2 + \left(-\frac{1}{2\sqrt{2}}\right)^2}}$$

$$< \frac{\sqrt{2}}{2}. \tag{11}$$

From (11) and the definition of $K$ and the definition of $\succcurlyeq$ given by (3), it implies

$$P_K(u) - P_K(v) \notin K.$$

That is

$$P_K(u) \not\succcurlyeq P_K(v). \tag{12}$$

By (10) and (12), it follows that the metric projection operator $P_K$ is not $\succcurlyeq$-increasing.

## 4. Isotone cones in Banach spaces and examples

In [9], Kong, Liu and Wu proved that, in a smooth, strictly convex and reflexive Banach lattice $(X, \succcurlyeq)$, where $\succcurlyeq$ is induced by a nonempty closed and convex cone $K \subseteq X$, if $K$ is orthogonal and subdual, then the metric projection operator $P_K$ satisfies

$$P_K(x) = x^+, \text{ for all } x \in X. \tag{13}$$

Since the operator $(\cdot)^+: X \to K$ is $\succcurlyeq$-increasing, it implies that the metric projection $P_K: X \to K$ is $\succcurlyeq$-increasing. In the same paper, the authors applied the $\succcurlyeq$-increasing property of $P_K$ to prove the existence of solution of some coupled best approximation problems in smooth and strictly convex reflexive Banach spaces.

In this section, we first introduce the concept of isotone cones in Banach spaces, then we provide some useful examples for applications in the related fields. Here the concept of "isotone cone" means that the metric projection operator onto a cone is isotone with respect to the partial order induced by this cone.

**Definition 4.1.** Let $X$ be a Banach space and $K$ a nonempty closed and convex cone in $X$ that induces a partial order $\succcurlyeq$ on $X$. $K$ is called a downward (upward) isotone cone in $X$, if

(a) $P_K(x) \neq \emptyset$, for every $x \in C$;
(b) $P_K$ is $\succcurlyeq$-increasing downward (upward).

$K$ is called an isotone cone in $X$, if it both of $\succcurlyeq$-increasing downward and upward. An isotone cone $K$ in $X$ is called a Chebyshev isotone cone, if $P_K$ is single-valued; that is,

(c) $K$ is a Chebyshev set;
(d) $P_K: X \to K$ is $\succcurlyeq$-increasing.

We list some examples below. These examples may be applied to some approximation problems related to the metric projection operators onto cones in Banach spaces.

**Example 4.2** (see [9]). In a Banach lattice, the operator $(\cdot)^+: X \to K$ is single-valued and order-increasing. From (13), it implies that every orthogonal and subdual cone in a smooth, strictly convex and reflexive Banach lattice is a Chebyshev isotone cone.

**Example 4.3.** Let $K$ be a nonempty closed and convex cone in a Banach space $X$. It is well-known that, for every $x \in X$, if $X$ is reflexive, then, $P_K(x) \neq \emptyset$; if $X$ is strictly convex, then, $P_K(x)$ contains at most one point. It implies that every isotone cone in a strictly convex and reflexive Banach space must be a Chebyshev isotone cone.

**Example 4.4.** Let $a < b$ and $C[a, b]$ be the Banach space of continuous functions on $[a, b]$ with absolute maximum norm $\|\cdot\|$. Let $C^+ = C^+[a, b]$ denote the positive cone in $C[a, b]$. That is

$$C^+ = \{x \in C[a, b]: x(t) \geq 0, \text{ for all } t \in [a, b]\}.$$

Then $C^+$ is a nonempty closed and convex cone in $C[a, b]$. Let $\succcurlyeq$ be the partial order on $C[a, b]$ induced by $C^+$. It satisfies, for $x, y \in C[a, b]$,

$$x \preccurlyeq y \text{ if and only if } x(t) \leq y(t), \text{ for all } t \in [a, b].$$

For arbitrary $x, y \in C[a, b]$, one can check that $x \vee y$ exists and it is calculated by

$$(x \vee y)(t) = \max\{x(t), y(t)\}, \text{ for all } t \in [a, b].$$

If $0 \preccurlyeq x \preccurlyeq y$ in $C[a, b]$, it is clear to have $\|x\| \leq \|y\|$. It implies that $(C[a, b], \succcurlyeq)$ is a Banach lattice (Hence, it is a normal partially ordered Banach space). We can check that $(C[a, b], \succcurlyeq)$ is not chain-complete ("chain-complete" means that every chain that is $\succcurlyeq$-bounded from above has its least $\succcurlyeq$-upper bound). It implies that $(C[a, b], \succcurlyeq)$ is not $\sigma$-Dedekind complete, so it is not Dedekind complete.

It is well-known that $C[a, b]$ is neither reflexive, nor strictly convex. The unreflexively property of $C[a, b]$ causes that, in general, the metric projection operator from $C[a, b]$ onto a nonempty closed and convex cone may be or may not be well defined. Fortunately, for this special positive cone $C^+$ in $C[a, b]$, we can show that $P_{C^+}: C[a, b] \to 2^{C^+}\setminus\{\emptyset\}$ is a well-defined set-valued mapping and we will deduce a formula for calculating the values of $P_{C^+}$. To this end, we need some notations: for $x \in C[a, b]$, we write

$$\zeta(x) = \{t \in [a, b]; x(t) \leq 0\}$$

and

$$\lambda(x) = \begin{cases} \max\{|x(t)|: t \in \zeta(x)\}, & \text{if } \zeta(x) \text{ is nonempty,} \\ 0, & \text{if } \zeta(x) \text{ is empty.} \end{cases} \quad (14)$$

It implies that

$$\lambda(x) = 0, \text{ for all } x \in C^+.$$

We write

$$\xi(x) = \begin{cases} \{t \in [a, b]: t \in \zeta(x) \text{ and } |x(t)| = \lambda(x)\}, & \text{if } \zeta(x) \text{ is nonempty,} \\ \text{empty}, & \text{if } \zeta(x) \text{ is empty.} \end{cases}$$

One can check that the metric projection operator $P_K$ satisfies that, for all $x \in C[a, b]$,

$$P_{C^+}(x) = \{z \in C^+: z(t) = 0, \text{ if } t \in \xi(x) \text{ and } |z(t) - x(t)| \leq \lambda(x), \text{ if } t \notin \xi(x)\}. \quad (15)$$

Hence

$$x^+ \in P_{C^+}(x), \text{ for all } x \in C[a, b], \quad (16)$$

where

$$x^+(t) = \begin{cases} x(t), & \text{if } t \notin \zeta(x), \\ 0, & \text{if } t \in \zeta(x). \end{cases} \quad (17)$$

It implies that, for all $x \in C[a, b], P_{C^+}(x) \neq \emptyset$. From (15), we have

$$\|z - x\| = \lambda(x), \text{ for } x \in C[a, b], \text{ and for every } z \in P_{C^+}(x).$$

Moreover, from (14) and (15), it implies that the values of $P_{C^+}$ have the following properties:

(a) $P_{C^+}(x) = \{x\}$, which is a singleton, for $x \in C^+$;

(b) $P_{C^+}(x)$ is an infinity set if $x \notin C^+$.

Next we show that $P_{C^+}$ is $\succcurlyeq$-increasing downward. For any given $x, y \in C[a, b]$ with $x \preccurlyeq y$, from the definition of $\succcurlyeq$ and (16), (17), it implies

$$\zeta(y) \subseteq \zeta(x), \quad \lambda(y) \le \lambda(x), \text{ and } x^+ \preccurlyeq y^+. \tag{18}$$

There is an $\succcurlyeq$-smallest element in $P_{C^+}(y)$, denoted by $\underline{y}$ that is

$$\underline{y}(t) = \min\{y^+(t) - \lambda(y), 0\}, \text{ for all } t \in [a, b].$$

Similarly, $P_{C^+}(x)$ has an $\succcurlyeq$-smallest element $\underline{x}$ and it is

$$\underline{x}(t) = \min\{x^+(t) - \lambda(x), 0\}, \text{ for all } t \in [a, b].$$

From (14), (17) and (18), we have $\underline{x}(t) \le \underline{y}(t)$ for all $t \in [a, b]$. It implies $\underline{x} \preccurlyeq \underline{y}$. Then, for any $z \in P_{C^+}(y)$, we pick element $\underline{x} \in P_{C^+}(y)$ which satisfies

$$z \succcurlyeq \underline{y} \succcurlyeq \underline{x}.$$

It implies that $P_{C^+}$ is $\succcurlyeq$-increasing downward (We can show that $P_{C^+}$ is not $\succcurlyeq$-increasing upward). It concludes that $C^+$ is a downward isotone cone in $C[a, b]$.

**Example 4.5.** Let $l_\infty$ be the Banach space of real sequences with absolute supremum norm. Let $l_\infty^+$ be the positive cone of $l_\infty$ (That is the set of all elements of $l_\infty$ with nonnegative components). Let $\succcurlyeq$ be the partial order on $l_\infty$ induced by $l_\infty^+$. For every $x = (x_1, x_2, \ldots)$, $y = (y_1, y_2, \ldots) \in l_\infty$, we have

$$x \preccurlyeq y \text{ if and only if } x_n \le y_n, \text{ for } n = 1, 2, \ldots . \tag{19}$$

For arbitrary $x = (x_1, x_2, \ldots)$, $y = (y_1, y_2, \ldots) \in l_\infty$, one can check that $x \vee y$ exists and it satisfies

$$(x \vee y)_n = \max\{x_n, y_n\}, \text{ for all } n = 1, 2, \ldots .$$

If $0 \preccurlyeq x \preccurlyeq y$ in $l_\infty$, it is clear to have $\|x\| \le \|y\|$. It implies that $(l_\infty, \succcurlyeq)$ is a Banach lattice. One can check that $(l_\infty, \succcurlyeq)$ is Dedekind complete, so it is chain-complete and $(l_\infty, \succcurlyeq)$ is a normal partially ordered Banach space (But it is not regular).

It is well-known that $l_\infty$ is neither reflexive, nor strictly convex. Similar to the last example, for this special positive cone $l_\infty^+$ in $l_\infty$, we can show that $P_{l_\infty^+} : l_\infty \to 2^{l_\infty^+} \setminus \{\emptyset\}$ is a well-defined set-valued mapping and we will deduce a formula for $P_{l_\infty^+}$. To this end, we use the same notations as in the last example: for $x = (x_1, x_2, \ldots) \in l_\infty$, we write

$$\zeta(x) = \{n \in N; x_n \le 0\}$$

and

$$\lambda(x) = \begin{cases} \max\{|x_n| : n \in \zeta(x)\}, & \text{if } \zeta(x) \text{ is nonempty}, \\ 0, & \text{if } \zeta(x) \text{ is empty}. \end{cases} \quad (20)$$

Then

$$\lambda(x) = 0, \text{ for all } x \in l_\infty^+.$$

We write

$$\xi(x) = \begin{cases} \{n \in N : n \in \zeta(x) \text{ and } |x_n| = \lambda(x)\}, & \text{if } \zeta(x) \text{ is nonempty}, \\ \text{empty}, & \text{if } \zeta(x) \text{ is empty}. \end{cases}$$

The metric projection operator $P_K$ satisfies that, for all $x = (x_1, x_2, \ldots) \in l_\infty$,

$$P_{l_\infty^+}(x) = \{z = (z_1, z_2, \ldots) \in K : z_n = 0, \text{ if } n \in \xi(x) \text{ and } |z_n - x_n| \le \lambda(x), \text{ if } n \notin \xi(x)\}. \quad (21)$$

Hence

$$x^+ \in P_{l_\infty^+}(x), \text{ for all } x \in l_\infty, \quad (22)$$

where

$$(x^+)_n = \begin{cases} x_n, & \text{if } n \notin \zeta(x), \\ 0, & \text{if } n \in \zeta(x). \end{cases}$$

It implies that, for all $x \in l_\infty$, $P_{l_\infty^+}(x) \ne \emptyset$. From (21), we have

$$\|z - x\| = \lambda(x), \text{ for } x \in l_\infty, \text{ and for every } z \in P_{l_\infty^+}(x).$$

Moreover, from (20) and (21), it implies that the values of $P_{l_\infty^+}$ have the following properties:

(a) $P_{l_\infty^+}(x) = \{x\}$, which is a singleton, for $x \in l_\infty^+$;
(b) $P_{l_\infty^+}(x)$ is an infinity set if $x \notin l_\infty^+$.

Next we show that $P_{l_\infty^+}$ is ≽-increasing downward. Take arbitrary $x = (x_1, x_2, \ldots)$, $y = (y_1, y_2, \ldots)$ $\in l_\infty$ satisfying $x \preccurlyeq y$, from (19) and (22), we have

$$\zeta(y) \subseteq \zeta(x), \quad \lambda(y) \le \lambda(x) \quad \text{and} \quad x^+ \preccurlyeq y^+. \quad (23)$$

In $P_{l_\infty^+}(y)$, there is ≽-smallest element denoted by $\underline{y}$ and it is

$$\underline{y}_n = \min\{(y^+)_n - \lambda(y), 0\}, \text{ for } n = 1, 2, \ldots. \quad (24)$$

$P_{l_\infty^+}(x)$ has the ≽-smallest element $\underline{x}$ and it is

$$\underline{x}_n = \min\{(y^+)_n - \lambda(x), 0\}, \text{ for } n = 1, 2, \ldots. \tag{25}$$

From (23)–(25), we have

$$\underline{x}_n \leq \underline{y}_n, \text{ for } n = 1, 2, \ldots.$$

It implies

$$\underline{x} \preccurlyeq \underline{y}.$$

Then, for any $z \in P_{l_\infty^+}(y)$, we pick the element $\underline{x} \in P_{l_\infty^+}(x)$ that satisfies

$$z \succcurlyeq \underline{y} \succcurlyeq \underline{x}.$$

It implies that $P_{l_\infty^+}$ is $\succcurlyeq$-increasing downward and $l_\infty^+$ is a downward isotone cone in $l_\infty$.

**Example 4.6.** Let $(\Sigma, \mu)$ be a $\sigma$-finite measure space and let $L_\infty = L_\infty(\Sigma, \mu)$. Let $L_\infty^+$ be the positive cone of $L_\infty$ (It is the set of all elements of $L_\infty$ with nonnegative values almost everywhere on $\Sigma$) that induces the partial order $\succcurlyeq$ on $L_\infty$, which satisfies that, for $x, y \in L_\infty$,

$$x \preccurlyeq y \text{ if and only if } x(t) \leq y(t), \text{ for almost all } t \in \Sigma.$$

For arbitrary $x, y \in L_\infty$, one can check that $x \vee y$ exists and it satisfies

$$(x \vee y)(t) = \max\{x(t), y(t)\}, \text{ for all } t \in \Sigma.$$

If $0 \preccurlyeq x \preccurlyeq y$ in $L_\infty$, it is clear to have $\|x\| \leq \|y\|$. It implies that $(L_\infty, \succcurlyeq)$ is a Banach lattice. One can check that $(L_\infty, \succcurlyeq)$ is not Dedekind complete. But it is $\sigma$-Dedekind complete, so it is chain-complete.

Similar to the last example, we can show that $P_{L_\infty^+}$ is an $\succcurlyeq$-increasing downward set-valued mapping. It concludes that $L_\infty^+$ is a downward isotone cone in $L_\infty$.

**Example 4.7.** Let $(\Sigma, \mu)$ be a $\sigma$-finite measure space. For $1 \leq p < \infty$, we consider the Banach space $L_p = L_p(\Sigma, \mu)$. Let $L_p^+$ be the positive cone of $L_p$ (it is the set of all elements of $L_p$ with nonnegative values almost everywhere on $\Sigma$). $L_p^+$ induces the partial order $\succcurlyeq$ on $L_p$, which satisfies that, for $x, y \in L_p$,

$$x \preccurlyeq y \text{ if and only if } x(t) \leq y(t), \text{ for almost all } t \in \Sigma.$$

For arbitrary $x, y \in L_\infty$, one can check that $x \vee y$ exists and it satisfies

$$(x \vee y)(t) = \max\{x(t), y(t)\}, \text{ for all } t \in \Sigma.$$

If $0 \preccurlyeq x \preccurlyeq y$ in $L_p$, it is clear to have $\|x\| \leq \|y\|$. It implies that $(l_p, \succcurlyeq)$ is a Banach lattice. One can check that $(L_p, \succcurlyeq)$ is not Dedekind complete. But it is $\sigma$-Dedekind complete, so it is chain-complete. One can show that the metric projection operator $P_{L_p^+}$ satisfies

$$P_{L_p^+}(x) = x^+, \text{ for all } x \in L_p,$$

where

$$x^+(t) = \begin{cases} x(t), & \text{if } x(t) \geq 0, \\ 0, & \text{if } x(t) < 0, \end{cases} \text{ for all } t \in \Sigma.$$

It implies that, $P_{L_p^+}: L_p \to L_p^+$ is a single-valued mapping. Since the operator $(\cdot)^+: L_p \to L_p^+$ is ≽-increasing. So is $P_{L_p^+}$. We conclude that $L_p^+$ is a Chebyshev isotone cone in $L_p$. Furthermore, one can check that, for $1 < p < \infty$, $L_p^+$ is an orthogonal and subdual cone in $L_p$. In a different view, it also implies that $L_p^+$ is a Chebyshev isotone cone in $L_p$, which is a smooth, strictly convex and reflexive Banach lattice.

**Example 4.8.** For $1 \leq p < \infty$, let $l_p$ be the Banach space of real sequences with absolute supremum norm. Let $l_p^+$ be the positive cone of $l_p$ (That is the set of all elements of $l_p$ with nonnegative components). Let ≽ be the partial order on $l_p$ induced by $l_p^+$. For every $x = (x_1, x_2, \ldots)$, $y = (y_1, y_2, \ldots) \in l_p$, we have

$$x \preccurlyeq y \text{ if and only if } x_n \leq y_n, \text{ for } n = 1, 2, \ldots.$$

For arbitrary $x = (x_1, x_2, \ldots)$, $y = (y_1, y_2, \ldots) \in l_p$, one can check that $x \vee y$ exists and it satisfies

$$(x \vee y)_n = \max\{x_n, y_n\}, \text{ for all } n = 1, 2, \ldots.$$

If $0 \preccurlyeq x \preccurlyeq y$ in $l_p$, it is clear to have $\|x\| \leq \|y\|$. It implies that $(l_p, \succcurlyeq)$ is a Banach lattice. One can check that $(l_p, \succcurlyeq)$ is Dedekind complete, so it is chain-complete ($(l_p, \succcurlyeq)$ is a normal partially ordered Banach space. It is also regular).

It is well-known that $l_1$ is neither reflexive, nor strictly convex. Similar to the last example, for this positive cone $l_p^+$ in $l_p$, we have that $P_{l_p^+}: l_p \to l_p^+$ is a well-defined single-valued mapping satisfying

$$P_{l_p^+}(x) = x^+, \text{ for all } x \in l_p,$$

where

$$(x^+)_n = \begin{cases} x_n, & \text{if } n \notin \zeta(x), \\ 0, & \text{if } n \in \zeta(x). \end{cases}$$

It concludes that $l_p^+$ is a Chebyshev isotone cone in $l_p$.

## 5. Fixed point theorems of order-increasing downward mappings on re-chain complete posets

In the early publications on fixed point theory in partially ordered sets (posets), if the considered mappings are set-valued mappings, then the mappings are usually required to be order-increasing upward (see [10–12]). In the examples in the previous section, we find that the metric projection operator $P_K$ from a Banach space $X$ onto an isotone cone $K \subseteq X$ is $\succcurlyeq$-increasing downward, where the partial order $\succcurlyeq$ on $X$ is induced by $K$. It leads us to study fixed point theorem on posets for order-increasing downward mappings. It is expected that the following theorems and their proofs at this aspect are analogous to those for order-increasing upward mappings. We recall some concepts and notations (See [10–12]).

A poset $(P, \succcurlyeq)$ is said to be re-inductive if every chain has a lower $\succcurlyeq$-bound in $P$. $(P, \succcurlyeq)$ is said to be re-chain-complete if every chain has its greatest lower $\succcurlyeq$-bound in $P$.

Let $A$ be a nonempty subset of a poset $(P, \succcurlyeq)$. $A$ is said to be *universally re-inductive* in $P$ whenever it satisfies that, for any given chain $\{x_\alpha\} \subseteq P$, if every element $x_\beta \in \{x_\alpha\}$ has a lower cover in $A$, then $\{x_\alpha\}$ has a lower $\succcurlyeq$-bound in $A$.

We list some results about universally re-inductive subsets of posets without proofs. For the details and proofs, the readers are referenced to [10–12].

**Lemma 5.1**. *Every re-inductive subset $A$ in a poset with a finite number of minimal elements is universally re-inductive.*

**Lemma 5.2.** *Every non-empty compact subset of a partially ordered Hausdorff topological space is both universally inductive and universally re-inductive.*

**Lemma 5.3.** *Let $X$ be a Banach space. Let $\succcurlyeq$ be a partial order on $X$ induced by a nonempty closed and convex cone $K$ in $X$. Let $\omega$ be the weak topology on $X$. Then the norm topology is natural with respect to $\succcurlyeq$ if and only if the weak topology $\omega$ is natural with respect to $\succcurlyeq$.*

**Lemma 5.4**. *Every non-empty bounded, closed and convex subset of a reflexive partially ordered Banach space is both universally inductive and universally re-inductive.*

Let $(X, \succcurlyeq^X)$, $(U, \succcurlyeq^U)$ be posets and let $F: X \to 2^U \setminus \{\emptyset\}$ be a set-valued mapping. $F$ is said to be order-increasing upward, if $x \preccurlyeq^X y$ in $X$, then, for any $z \in F(x)$, there is a $w \in F(y)$ such that $z \preccurlyeq^U w$. $F$ is said to be order-increasing downward, if $x \preccurlyeq^X y$ in $X$, then, for any $w \in F(y)$, there is a $z \in F(x)$ such that $z \preccurlyeq^U w$. If $F$ is both order-increasing upward and downward, then $F$ is said to be order-increasing.

In [10–12], several fixed point theorems for order-increasing upward mappings were proved. In this section, we prove a fixed point theorem for order-increasing downward mappings.

**Theorem 5.1**. *Let $(P, \succcurlyeq)$ be a re-chain-complete poset and let $F: P \to 2^P \setminus \{\emptyset\}$ be a set-valued mapping satisfying the following three conditions*:

  (a). *$F$ is order-increasing downward*;
  (b). *$(F(x), \succcurlyeq)$ is universally re-inductive, for every $x \in P$*;
  (c). *There is an element $y_*$ in $P$ and $v_* \in F(y_*)$ with $v_* \preccurlyeq y_*$.*
*Then*

1. $(\mathcal{F}(F), \succcurlyeq)$ is a nonempty re-inductive poset;
2. $(\mathcal{F}(F) \cap (y_*], \succcurlyeq)$ is a nonempty re-inductive poset.

*Consequently, we have*

(i)   *F has an $\succcurlyeq$-minimal fixed point;*
(ii)  *F has an $\succcurlyeq$-minimal fixed point $x^*$ with $x^* \preccurlyeq y_*$.*

*Proof.* Define a subset $B \subseteq P$ by

$$B = \{z \in P: \text{there is } v \in F(z) \text{ such that } v \preccurlyeq z\}.$$

Since the point $y_* \in B$, it implies that $B \neq \emptyset$. Then we prove that $(B, \succcurlyeq)$ is re-chain-complete (it follows that $(B, \succcurlyeq)$ is re-inductive). Take an arbitrary chain $\{x_\alpha\} \subseteq B$. Since $(P, \succcurlyeq)$ is re-chain-complete, $\wedge\{x_\alpha\}$ exists in $P$. We write $x = \wedge\{x_\alpha\}$. Then we show that $x$ is a lower $\succcurlyeq$-bound of $\{x_\alpha\}$ in $B$. For every $\alpha$, from $x_\alpha \in B$, there is $v_\alpha \in F(x_\alpha)$ such that $v_\alpha \preccurlyeq x_\alpha$. From $x \preccurlyeq x_\alpha$ and condition (a), the order-increasing downward property of $F$, for $v_\alpha \in F(x_\alpha)$, there is $u_\alpha \in F(x)$ such that $u_\alpha \preccurlyeq v_\alpha$. It follows that

$$u_\alpha \preccurlyeq v_\alpha \preccurlyeq x_\alpha, \text{ for all } \alpha.$$

Since $F(x)$ is universally re-inductive, there is $u \in F(x)$ such that

$$u \preccurlyeq x_\alpha, \text{ for all } \alpha.$$

Hence $u$ is a lower $\succcurlyeq$-bound of $\{x_\alpha\}$ in $P$. On the other hand, from $u \in F(x)$ and $u \preccurlyeq x$ and the order-increasing downward property of $F$, there is $v \in F(u)$ such that $v \preccurlyeq u$. It implies that $u \in B$. Hence $u$ is a lower $\succcurlyeq$-bound of $\{x_\alpha\}$ in $B$. Since $x = \wedge\{x_\alpha\}$, it implies $u \preccurlyeq x$; and therefore, $x \in B$. It shows that

$$(B, \succcurlyeq) \text{ is re-chain-complete. So it is re-inductive.} \tag{26}$$

Take any $\succcurlyeq$-minimal point $y \in B$. From (26), there is $w \in F(y)$ such that $w \preccurlyeq y$. From the above proof, we must have $w \in B$. Since $w \preccurlyeq y$ and $y$ is an $\succcurlyeq$-minimal point in $B$, it implies $y = w \in F(y)$. Hence $y$ is a fixed point of $F$. That is, $\mathcal{F}(F) \neq \emptyset$. Meanwhile we obtain that

$$\text{every } \succcurlyeq\text{-minimal point in } B \text{ is a fixed point of } F. \tag{27}$$

Take an arbitrary chain $\{z_\beta\}$ in $\mathcal{F}(F)$. Since $\mathcal{F}(F) \subseteq B$, from (26), $\wedge\{z_\beta\}$ exists in $B$. We write $z = \wedge\{z_\beta\} \in B$. By (26) again, there is an $\succcurlyeq$-minimal point $a \in B$ such that $a \preccurlyeq z$. It implies that $a$ is a lower $\succcurlyeq$-bound of $\{z_\beta\}$. From (27), $a \in \mathcal{F}(F)$. It concludes that $(\mathcal{F}(F), \succcurlyeq)$ is re-inductive. Hence it proves part 1. To prove part 2, we define a restricted set-valued mapping $F^*$ on $(y_*]$ by

$$F^*(x) = F(x) \cap (y_*], \text{ for } x \in (y_*]. \tag{28}$$

It is clear that $((y_*], \succcurlyeq)$ is a re-chain-complete poset (a sub-poset of a re-chain-complete poset $(P, \succcurlyeq)$). From (28), it follows that

$$\mathcal{F}(F^*) = \mathcal{F}(F) \cap (y_*] \subseteq \mathcal{F}(F). \tag{29}$$

From condition (c), $v_* \in F(y_*)$ and $v_* \preccurlyeq y_*$, it implies $v_* \in F^*(y_*)$ and $F^*(y_*) \neq \emptyset$. Then we show that, for every for $x \in (y_*]$, $F^*(x) \neq \emptyset$. From condition (a) again, the $\succcurlyeq$-increasing downward property of $F$, for $x \in (y_*]$ and $v_* \in F(y_*)$, there is $u \in F(x)$ such that $u \preccurlyeq v_* \preccurlyeq y_*$. It implies $u \in F^*(x)$ and $F^*(x) \neq \emptyset$.

Next we prove that $F^*$ is $\succcurlyeq$-increasing downward on $(y_*]$. For any $x, y \in (y_*]$ with $x \preccurlyeq y$, take any $v \in F^*(y) \subseteq F(y)$, from the $\succcurlyeq$-increasing downward property of $F$, there is $w \in F(x)$ such that $w \preccurlyeq v \in F^*(y) = F(y) \cap (y_*]$. It implies that $w \in (y_*]$; and therefore, $w \in F^*(x)$. So $F^*$ is $\succcurlyeq$-increasing downward on $(y_*]$.

From condition (b), for every $x \in P$, $(F(x), \succcurlyeq)$ is universally re-inductive. It implies that, for every $x \in (y_*]$, $F^*(x) = F(x) \cap (y_*]$ must be universally re-inductive.

Then by using the restricted mapping $F^*$, part 2 follows immediately from the conclusion part 1 and (29). □

By using Lemmas 5.2, 5.3, and 5.4, we obtain the following corollary of Theorem 5.1.

**Corollary 5.1**. *Let $(X, \succcurlyeq)$ be a partially ordered reflexive Banach space, in which the partial order $\succcurlyeq$ is induced by a nonempty closed and convex cone $K$, and $D$ a nonempty bounded closed convex subset of $X$. Let $F: D \rightarrow 2^D \setminus \{\emptyset\}$ be a set-valued mapping satisfying the following conditions*:

(a). *$F$ is $\succcurlyeq$-increasing downward*;
(b). *$F(x)$ is a nonempty closed and convex subset of $D$, for every $x \in D$*;
(c). *There is an element $y_*$ in $P$ and $v_* \in F(y_*)$ with $v_* \preccurlyeq y_*$.*

*Then*

1. *$(\mathcal{F}(F), \succcurlyeq)$ is a nonempty re-inductive poset*;
2. *$(\mathcal{F}(F) \cap (y_*], \succcurlyeq)$ is a nonempty re-inductive poset.*

*Consequently, we have*

(i)   *$F$ has an $\succcurlyeq$-minimal fixed point*;
(ii)  *$F$ has an $\succcurlyeq$-minimal fixed point $x_*$ with $x_* \preccurlyeq y_*$.*

## 6. Best approximation theorems for operators without continuity conditions

In this section, we apply the fixed point theorem given in the previous section and the concepts of isotone cones in partially ordered Banach spaces to solve some best approximation problems for some operators, which are not required to have any continuity property.

In Theorem 6.1 in the previous section, the underlying poset $(P, \succcurlyeq)$ is required to be re-chain-complete. When we apply this theorem to partially ordered Banach spaces, we need some concepts related to chain-complete and re-chain-complete.

We say that a partially ordered Banach space is chain-complete, if every chain with an upper bound has the smallest upper bound and it is re-chain-complete if every chain with a lower bound has the greatest lower bound.

We list some examples of chain-complete and re-chain-complete partially ordered Banach spaces. For the details and proofs, the readers are referenced to [10–12].

**Lemma 6.1**. *Every regular partially ordered Banach space is both chain-complete and re-chain-complete.*

**Lemma 6.2.** *Every normal partially ordered reflexive Banach space (it is also regular) is both chain-complete and re-chain-complete.*

**Lemma 6.3.** *Every reflexive Banach lattice (it is both normal and regular) is both chain-complete and re-chain-complete.*

Let $X$ be a Banach space and $K$ a nonempty closed and convex cone in $X$. Let $f: K \to X$ be a mapping. For any given point $x \in K$, an element $y \in K$ is called a best approximation point for $f(x)$ on $K$ if

$$\|f(x) - y\| = \inf\{\|f(x) - z\|: z \in K\}.$$

The set of best approximation points for $f(x)$ on $K$ is

$$P_K(f(x)) = \{y \in K: \|f(x) - y\| = \inf\{\|f(x) - z\|: z \in K\}\}.$$

For a point $x_0 \in K$, if $x_0 \in P_K(f(x_0))$, then $x_0$ is called a best approximation point for $f$ on $K$. It satisfies

$$\|f(x_0) - x_0\| = \inf\{\|f(x_0) - z\|: z \in K\}. \tag{30}$$

Let $\mathcal{B}(f)$ denotes the set of best approximation points for $f$ on $K$.

**Theorem 6.1**. *Let $K$ be a downward isotone cone in a Banach space $X$ that induces a partial order $\succcurlyeq$ on $X$. Let $f: K \to X$ be a mapping. Suppose that the following conditions are satisfied*:

a) $(X, \succcurlyeq)$ *is re-chain-complete*;
b) $(P_K(x), \succcurlyeq)$ *is universally re-inductive, for every $x \in X$*;
c) $f$ *is $\succcurlyeq$-increasing*;
d) *There is $y_* \in K$ such that $f(y_*) \preccurlyeq y_*$.*

*Then $f$ has a best approximation point on $K$. Moreover,*

$$(\mathcal{B}(f) \cap (y_*], \succcurlyeq) \text{ is a nonempty re-inductive poset.}$$

*Consequently, $f$ has an $\succcurlyeq$-minimal best approximation point $x^*$ on $K$ satisfying $x^* \preccurlyeq y_*$.*

*Proof.* Since $K$ is a downward isotone cone in $X$, $P_K: X \to 2^K \setminus \{\emptyset\}$ is an $\succcurlyeq$-increasing downward set-valued mapping. From condition d), $f(y_*) \preccurlyeq y_*$, and noticing that $P_K(y_*) = y_*$, it implies that there is $v_* \in P_K(f(y_*))$ such that

$$v_* \preccurlyeq y_*. \tag{31}$$

Let $K(y_*) = K \cap (y_*]$. From condition a) in this theorem, $(K \cap (y_*], \succcurlyeq)$ is re-chain-complete. Define a set-valued mapping $F: K(y_*) \to 2^{K(y_*)} \setminus \{\emptyset\}$ by

$$F(x) = P_K(f(x)) \cap (y_*] = \{y \in K(y_*): \|f(x) - y\| = \inf\{\|f(x) - z\|: z \in K\}\}, \text{ for } x \in K(y_*). \quad (32)$$

By (31), we have

$$v_* \in F(y_*) \text{ such that } v_* \preccurlyeq y_*.$$

Hence $F$ satisfies condition (c) in Theorem 5.1. Similar to the proof of part 2 of Theorem 5.1, we can show that $F$ is a well-defined set-valued mapping on $K(y_*)$ with nonempty values and it is $\succcurlyeq$-increasing downward. Hence $F$ satisfies condition (a) in Theorem 5.1.

From condition b) in this theorem, for every $x \in X$, $(P_K(x), \succcurlyeq)$ is universally re-inductive. It implies that, for every $x \in X$, $P_K(x) \cap (y_*]$ must be universally re-inductive. Hence, for every $x \in K$, $F(x) = P_K(f(x)) \cap (y_*]$ must be universally re-inductive. So $F$ satisfies condition (b) in Theorem 5.1. By Theorem 5.1, $F$ has a fixed point. That is, $\mathcal{F}(F) \neq \emptyset$.

Next we show that $\mathcal{B}(f) \cap (y_*] = \mathcal{F}(F)$. For an arbitrary $x^* \in \mathcal{F}(F)$, from (32), we have

$$x^* \in F(x^*) = P_K(f(x^*)) \cap (y_*] = \{y \in K(y_*): \|f(x^*) - y\| = \inf\{\|f(x^*) - z\|: z \in K\}\}.$$

That is

$$\|f(x^*) - x^*\| = \inf\{\|f(x^*) - z\|: z \in K\}.$$

From (30), it implies $x^* \in \mathcal{B}(f) \cap (y_*]$. We obtain

$$\mathcal{F}(F) \subseteq \mathcal{B}(f) \cap (y_*]. \quad (33)$$

It can be similarly shown that

$$\mathcal{B}(f) \cap (y_*] \subseteq \mathcal{F}(F). \quad (34)$$

(33) and (34) imply

$$\mathcal{F}(F) = \mathcal{B}(f) \cap (y_*]. \quad (35)$$

Then, by (35), this theorem follows from the conclusion parts 1 and ii) of Theorem 5.1 immediately. □

We need the following lemma to prove a corollary of Theorem 6.1.

**Lemma 6.4.** *Let $K$ be a nonempty closed and convex cone in a Banach space $X$. For every $x \in X$, if $P_K(x) \neq \emptyset$, then $P_K(x)$ is a closed, bounded and convex subset of $K$.*

*Proof.* The proof of this lemma is straight forward and it is omitted here.

**Corollary 6.1.** *Let $K$ be a nonempty closed and convex cone in a reflexive Banach space $X$ that induces a partial order $\succcurlyeq$. Let $f: K \to X$ be a mapping. Suppose that the following conditions are satisfied*:

  a) $(X, \succcurlyeq)$ *is re-chain-complete*;

b) $f$ is $\succcurlyeq$-increasing;
c) There is $y_* \in K$ such that $f(y_*) \preccurlyeq y_*$;
d) $P_K$ is $\succcurlyeq$-increasing downward.

Then $f$ has a best approximation point on $K$. Moreover,

$$(\mathcal{B}(f) \cap (y_*], \succcurlyeq) \text{ is a nonempty re-inductive poset.}$$

Consequently, $f$ has an $\succcurlyeq$-minimal best approximation point $x_*$ on $K$ satisfying $x_* \preccurlyeq y_*$.

*Proof.* Since $X$ is reflexive, for every $x \in X$, $P_K(x) \neq \emptyset$. Then from condition d) in this theorem, it implies that $K$ is a downward isotone cone in $X$.

From Lemma 5.1, for every $x \in X$, $P_K(f(x))$ is a nonempty bounded closed and convex subset of $K$. From Lemma 5.4, it implies that $(F(x), \succcurlyeq)$ is a universally re-inductive subset of $K$. So $F$ satisfies all conditions in Theorem 6.1 and this corollary is proved. □

Next we prove a theorem for the existence of best approximations for upward isotone cones in Banach spaces.

**Theorem 6.2.** *Let $K$ be a upward isotone cone in a Banach space $X$ that induces a partial order $\succcurlyeq$ on $X$. Let $f: K \to X$ be a mapping. Suppose that the following conditions are satisfied:*

a) $(X, \succcurlyeq)$ is chain-complete;
b) $(P_K(x), \succcurlyeq)$ is universally inductive, for every $x \in X$;
c) $f$ is $\succcurlyeq$-increasing;
d) There is $y_* \in K$ such that $P_K(f(y_*)) \subseteq (y_*]$.

*Then $f$ has a best approximation point on $K$. Moreover,*

$$(\mathcal{B}(f) \cap (y_*], \succcurlyeq) \text{ is a nonempty inductive poset.}$$

*Proof.* Similar to the proof of Theorem 6.1, let $K(y_*) = K \cap (y_*]$. From condition a) in this theorem, $(K(y_*), \succcurlyeq)$ is chain-complete. The assumption that $P_K$ is $\succcurlyeq$-increasing upward and condition d) in this theorem imply that,

$$P_K(f(x)) \subseteq K(y_*), \text{ for } x \in K(y_*).$$

It implies

$$P_K(f(x)) \cap (y_*] = P_K(f(x)), \text{ for } x \in K(y_*). \quad (36)$$

Define a set-valued mapping $F: K(y_*) \to 2^{K(y^*)} \setminus \{\emptyset\}$ by

$$F(x) = P_K(f(x)) = \{y \in K: \|f(x) - y\| = \inf\{\|f(x) - z\|: z \in K\}\}, \text{ for } x \in K(y_*).$$

By (35), we have $F(x) \subseteq K(y_*)$, and $F: K(y_*) \to 2^{K(y^*)} \setminus \{\emptyset\}$ is well-defined.

From condition b), for every $x \in X$, $(P_K(x), \succcurlyeq)$ is universally inductive. From (36), it implies that, for every $x \in X$, $(F(x), \succcurlyeq)$ must be universally inductive. From $0 \preccurlyeq v$, for every $v \in F(0)$, it implies that F satisfies all conditions A3 in Theorem 3.1 in [12] on $(K(y_*), \succcurlyeq)$. Then this theorem follows immediately from the conclusion part 1 in Theorem 3.1 in [12]. □

Consequently, we have

**Corollary 6.2**. *Let K be a nonempty closed and convex cone in a reflexive Banach space X that induces a partial order $\succcurlyeq$. Let $f: K \to X$ be a mapping. Suppose that the following conditions are satisfied*:

   a) $(X, \succcurlyeq)$ *is chain-complete*;
   b) *f is $\succcurlyeq$-increasing*;
   c) *There is $y_* \in K$ such that $P_K(f(y_*)) \subseteq (y_*]$*;
   d) *$P_K$ is $\succcurlyeq$-increasing upward*.

*Then f has a best approximation point on K. Moreover,*

$$(\mathcal{B}(f) \cap (y_*], \succcurlyeq) \text{ is a nonempty inductive poset}.$$

## 7. Orthogonal and subdual cones in Banach lattices

As stated in Example 4.2 in section 4, in [9], Kong, Liu and Wu introduced the concepts of orthogonal and subdual cones in Banach lattices. They applied these concepts for solving some best approximation problems in smooth, strictly convex and reflexive Banach lattices. In this section, we extend the results about best approximation problems in [9] from smooth, strictly convex and reflexive Banach lattices to general Banach lattices.

We generalize the concepts of orthogonal and subdual cones to more general Banach lattices.

Let $(X, \succcurlyeq)$ be a Banach lattice, where $\succcurlyeq$ is induced by a nonempty closed and convex cone $K \subseteq X$. K is said to be orthogonal if, for any $x, y \in X$,

$$x \wedge y = 0 \text{ implies } \langle jx, y \rangle = 0, \text{ for every } jx \in Jx.$$

K is said to be subdual if, for any $x \in X$, and for every $jx \in Jx$, we have

$$\langle jx, z \rangle \geq 0, \text{ for every } z \in K.$$

**Lemma 7.1**. *Let $(X, \succcurlyeq)$ be a Banach lattice, where $\succcurlyeq$ is induced by a nonempty closed and convex cone $K \subseteq X$. If K is orthogonal and subdual, then, for any $x \in X$, we have $x^+ \in P_K(x)$.*

*Proof.* The proof is similar to the first part of the proof of Theorem 3.1 in [9]. From $x^- = x^+ - x$ and $x^- \wedge x^+ = 0$, for any $j(x^-) \in J(x^-) = J(x^+ - x)$, we have

$$\langle z - x^+, j(x^-) \rangle = \langle z, j(x^-) \rangle \geq 0, \text{ for all } z \in K.$$

From Lemma 2.1, it implies that $x^+ \in P_K(x)$. □

**Theorem 7.1**. *Let $(X, \succcurlyeq)$ be a Banach lattice, where $\succcurlyeq$ is induced by an orthogonal and subdual cone $K \subseteq X$. If K is orthogonal and subdual, then, for any $x \in X$, we have $x^+ \in P_K(x)$. Let $f: K \to X$ be a mapping. Suppose that the following conditions are satisfied*:

a) $(X, \succcurlyeq)$ is re-chain-complete;
b) $f$ is $\succcurlyeq$-increasing;
c) There is $y_* \in K$ such that $f(y_*) \preccurlyeq y_*$.

*Then $f$ has a best approximation point on $K$.*

*Proof.* Similar to the proof of Theorem 6.1, let $K(y_*) = K \cap (y_*] = [0, y_*]$. From condition a) in this theorem, $(K \cap (y_*], \succcurlyeq)$ is re-chain-complete. Define a single-valued mapping $F: K(y_*) \to K(y_*)$ by

$$F(x) = f(x)^+, \text{ for } x \in K(y_*).$$

From condition c), $f(y_*) \preccurlyeq y_*$, it implies

$$F(y_*) = f(y_*)^+ \preccurlyeq y_*. \tag{37}$$

From condition b), for any $x \in K(y_*)$, it implies that $f(x) \preccurlyeq f(y_*)$. Then we get

$$0 \preccurlyeq f(x)^+ \preccurlyeq f(y_*)^+ \preccurlyeq y_*, \text{ for } x \in K(y_*).$$

It implies that the mapping $F: K(y_*) \to K(y_*)$ is well defined. Since $(\cdot)^+: X \to K$ is $\succcurlyeq$-increasing, it yields that $F$ is a composite of two $\succcurlyeq$-increasing mappings. So $F$ is $\succcurlyeq$-increasing. By (37), $F$ satisfies all conditions of Theorem 5.1. Hence $F$ has a fixed point $x_*$. That is,

$$x_* = F(x_*) = f(x_*)^+.$$

From Lemma 7.1, we have

$$x_* = f(x_*)^+ \in P_K(f(x_*)).$$

It implies that $x_*$ is a best approximation point of $f$ on $K$. □

## 8. Applications to variational inequalities in Banach spaces

In this section, we consider some variational inequalities in which the underlying spaces are isotone cones in Banach spaces. As applications of the best approximation theorems provided in the previous section, we prove the solvability of some variational inequalities in Banach spaces.

**Definition 8.1**. Let $X$ be a Banach space and $C$ a nonempty closed and convex subset of $X$. Let $f: C \to X$ be a mapping. The variational inequality problem associated with $f$, $C$, denoted by VI($f$, $C$), is formulated as follows:

$$\text{find } x_* \in C \text{ and } j(x_* - f(x_*)) \in J(x_* - f(x_*))$$

such that

$$\langle z - x_*, j(x_* - f(x_*)) \rangle \geq 0, \text{ for all } z \in C.$$

The set of solutions of VI($f$, $C$) is denoted by $\mathcal{S}(f, C)$.

**Theorem 8.1**. *Let $K$ be a downward isotone cone in a Banach space $X$ that induces a partial order $\succcurlyeq$ on $X$. Let $f: K \to X$ be a mapping. Suppose $f$ and $P_K$ satisfy conditions a)–d) in Theorem 6.2. Then the problem VI($f$, $K$) has a solution. Moreover,*

$(\mathcal{S}(f, K) \cap (y_*], \succcurlyeq)$ *is a nonempty re-inductive poset.*

*Consequently, the problem* VI$(f, K)$ *has an $\succcurlyeq$-minimal solution $x_*$ satisfying $x_* \preccurlyeq y_*$.*

*Proof.* As the definition (32) in the proof of Theorem 6.1, we define the set-valued mapping $F$ on $K \cap (y_*]$ by

$$F(x) = P_K(f(x)) \cap (y_*], \text{ for } x \in K \cap (y_*].$$

From Lemma 2.1, for any $x \in X$, $x' \in P_K(x)$ if and only if, there is $j(x' - x) \in J(x' - x)$ such that

$$\langle z - x', j(x' - x) \rangle \geq 0, \text{ for all } z \in K.$$

By (36), for $x_* \in K$, $x_* \in \mathcal{S}(f, C) \cap (y_*]$, if and only if $x_* \in P_K(f(x_*)) = F(x_*)$; that is, $x_* \in \mathcal{F}(F)$. From (35) in the proof of Theorem 6.1, it implies

$$\mathcal{S}(f, C) \cap (y_*] = \mathcal{F}(F) = \mathcal{B}(f) \cap (y_*].$$

Then this theorem follows immediately from Theorem 6.1. □

As a consequence of Theorem 8.1 and Corollary 6.1, we have

**Corollary 8.1**. *Let $K$ be a nonempty closed and convex cone in a reflexive Banach space $X$ that induces a partial order $\succcurlyeq$. Let $f: K \to X$ be a mapping. Suppose that $f$ and $P_K$ satisfy conditions a)–d) in Corollary 6.1. Then the problem* VI$(f, K)$ *has a solution. Moreover,*

$(\mathcal{S}(f, K) \cap (y_*], \succcurlyeq)$ *is a nonempty re-inductive poset.*

*Consequently, the problem* VI$(f, K)$ *has an $\succcurlyeq$-minimal solution $x_*$ satisfying $x_* \preccurlyeq y_*$.*

Similar to Theorem 8.1 and Corollary 8.1, we can apply Theorem 6.2 and Corollary 6.2 to study the problem VI$(f, K)$ for upward isotone cones in Banach spaces.